# Determinants of Population Growth in Rajasthan: An Analysis

Prof. V.V. Singh[1]
Dr. Alka Mittal[2]
Dr. Neetish Sharma[3]
Prof. F. Smarandache[4]

## Abstract

Rajasthan is the biggest State of India and is currently in the second phase of demographic transition and is moving towards the third phase of demographic transition with very slow pace. However, state's population will continue to grow for a time period. Rajasthan's performance in the social and economic sector has been poor in past. The poor performance is the outcome of poverty, illiteracy and poor development, which co-exist and reinforce each other. There are many demographic and socio-economic factors responsible for population growth. This paper attempts to identify the demographic and socio-economic variables, which are responsible for population growth in Rajasthan with the help of multivariate analysis.

## 1. Introduction:

*Prof. Stephan Hawking (Cambridge University) was on Larry King Live. Larry King called him the "most intelligent person in the world". King asked some very key questions, one of them was: "what worries you the most?" Hawking said, "My biggest worry is population growth, and if it continues at the current rate, we will be standing shoulder to shoulder in 2600. Something has to happen, and I don't want it to be a disaster".*

The importance of population studies in India has been recognized since very ancient times. The 'Arthashastra' of Kautilya gives a detailed description of how to conduct a population, economic and agricultural census. During the reign of Akbar, Abul Fazal compiled the Ain-E-Akbari containing comprehensive data on population, industry, wealth and characteristics of population. During the British period, system of decennial census started with the first census in 1872.

The population growth of a region and its economic development are closely linked. India has been a victim of population growth. Although the country has achieved progress in the economic field, the population growth has wrinkled the growth potential. The need to check the population growth was realized by a section of the intellectual elite even before independence. Birth control was accepted by this group but implementation was restricted to the westernized minority in the cities. When the country attained independence and planning was launched, population control became one of the important items on the agenda of development. The draft outline of the First Five Year Plan said, "the increasing pressure of population on natural resources retards economic progress and limits seriously the rate of extension of social services, so essential to civilized existence."

India was one of the pioneers in health service planning with a focus on primary health care. Improvement in the health status of the population has been one of the major thrust areas for the social development programs of the country in the five year plans.  India is a signatory to the Alma

---

[1]  Prof. V.V. Singh isProfressor &  Head, Deptt. of Economics, University of Rajasthan, Jaipur.
[2]  Dr. Alka Mittal is Lecturer of Economics in S.S. Jain Subodh PG Mahila Mahavidhyalaya, Jaipur.
[3]  Dr. Neetish Sharma is Assistant Director (Plan-Finance) in Planning Department, Govt. of Rajasthan, Jaipur.
[4]  Prof. F. Smarandache, University of New Mexico, Gallup, NM 87301, USA

*Views expressed here are personal views of above, not concerned to their organizations.*

Ata Declaration (1978) whereby a commitment was made to achieve 'Health for All' by 2000 AD. We are in the end of the first decade of the 21st century but still have to go a long way to achieve this target. Rajasthan is lagging behind the all India average in the key parameters i.e. CBR, CDR, IMR, TFR & CPR. The state has made consistent efforts to improve quality of its people through improvement in coverage & quality of health care and implementation of disease control programs but the goals remain elusive due to high levels of fertility and mortality. According to the Report of the Technical Group on Population Projections, India will achieve the target of TFR = 2.1 (Net Reproduction Rate = 1) in 2026. Kerala & Tamilnadu had already achieved it in 1988 & 1993 respectively but Rajasthan will achieve it in 2048 & Uttar Pradesh in 2100.

Rajasthan is the largest state of the country with its area of 342239 sq. kms., which constitutes about 10.41% of the total area of the country. According to 2001 census, its population is 56.51 million. It consist 5.5 % population and ranks eighth in the country. In 1901, population of Rajasthan was 10.29 millions. In 1951, it reached to 15.97 millions with its slow growth during 1901-1951. Figure 1 shows that it increased rapidly after 1951. It reached to 34.26 million in 1981 and to 56.51 million in 2001. It has multiplied 5.5 times since 1901 and 3.5 times since 1951. Figure 2 shows decennial growth in population of the state. Before 1951, it increased by less than 20% growth per decade. In 1971-81, it shows the maximum growth rate of 32.97%. In 1981-91, it decreased by 4.53 percentage points and grew by 28.44%. The decade of 1991-2001 shows growth of 28.41%.

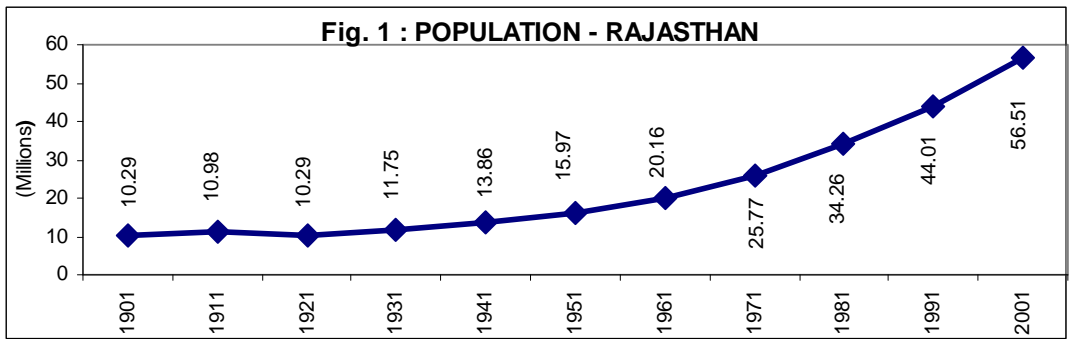

Source: Government of India, Registrar General, India, see the website www.censusindia.net

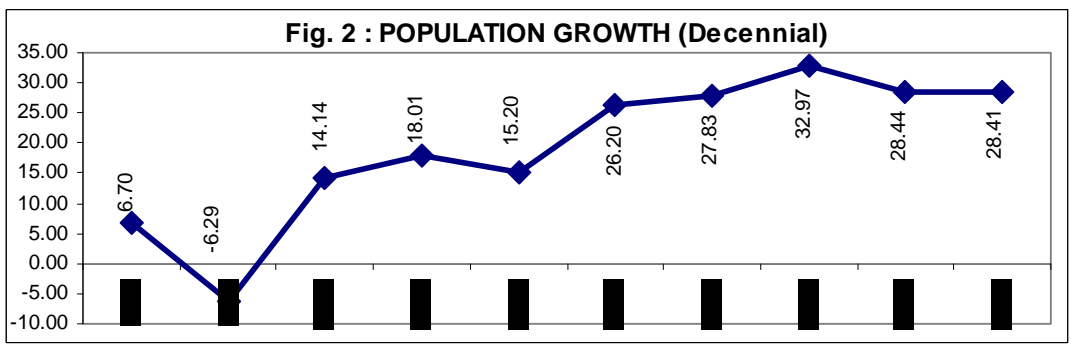

Source: Government of India, Registrar General, India, see the website www.censusindia.net

The rapid population growth in a already populated state like Rajasthan could lead to many problems i.e. pressure on land, environmental deterioration, fragmentation of land holding, shrinking forests, rising temperatures, pressure on health & educational infrastructure, on availability of food grains & on employment. Figure 3 shows the decennial growth of district-wise population during 1991-2001. Jaisalmer shows the maximum growth of 47.45% followed by Bikaner (38.18%), Barmer (36.83%), Jaipur (35.10%) and Jodhpur (33.77%). Rajasamand shows minimum growth of 19.88% followed by Jhunjhunu (20.90%), Chittorgarh (21.46%), Pali (22.39%) and Jhalawar (23.34%).

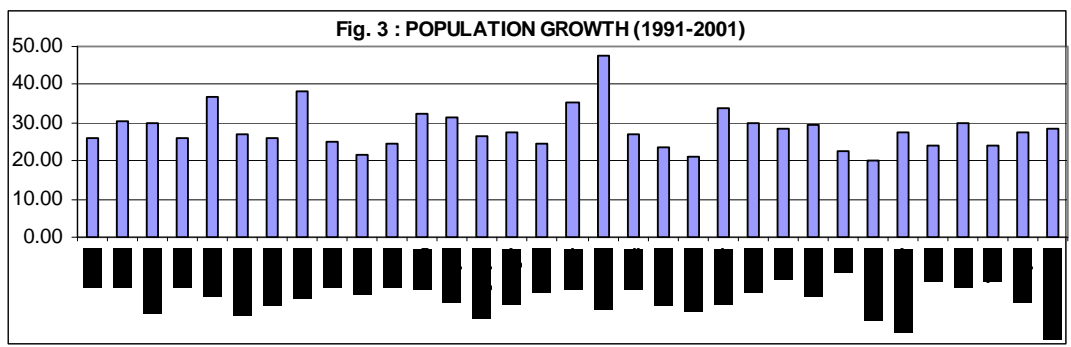

Source: Government of India, Registrar General, India, see the website www.censusindia.net

Rajasthan is currently in the second phase and is moving towards the third phase of demographic transition with very slow pace. The changes in the population growth rates in Rajasthan have been relatively slow, but the change has been steady and sustained. We are aware of the need for birth control, but too many remain ignorant of contraception methods or are unwilling to discuss them. There is considerable pressure to produce a son. However, the state's population will continue to grow for a time period.

Rajasthan is the second state in the country to formulate and adopt its own Population Policy in January 2000. State Population Policy[5] has envisaged strategies for population stabilization and improving health conditions of people specially women and children. The policy document has clearly presented role and responsibilities of different departments actively contributing in implementation of population policy. Family Welfare Program was linked with other sectors and demands intervention and efficient policies in these sectors so that changes can be brought in the social, economic, cultural & political environment. The State Population Policy envisages time bound objectives as mentioned in table 1:

**Table 1: Objectives of Population Policy of Rajasthan**

| Indicators | 1997 | 2001 | 2004 | 2007 | 2011 | 2013 | 2016 |
|---|---|---|---|---|---|---|---|
| **Total Fertility Rate** | 4.11 | 3.74 | 3.41 | 3.09 | 2.65 | 2.43 | 2.10 |
| **Birth Rate** | 32.1 | 29.2 | 27.5 | 25.6 | 22.6 | 20.9 | 18.4 |
| **Contraceptive Prevalence Rate** | 38.5 | 42.2 | 48.2 | 52.7 | 58.8 | 61.8 | 68.0 |
| **Death Rate** | 8.9 | 8.7 | 8.4 | 7.9 | 7.5 | 7.2 | 7.0 |
| **Infant Mortality Rate** | 85.0 | 77.4 | 72.7 | 68.1 | 62.2 | 60.1 | 56.8 |

Rajasthan's performance in the social and economic sector has been poor in past. The poor performance is the outcome of poverty, illiteracy and poor development which co-exist and reinforce each other. State Government has taken energetic steps in last few years to assess and fully meet the unmet needs for maternal & child health care and contraception through improvement in availability and access to family welfare services but still remains a long path. The progress in these indicators would determine the year and size of the population at which the state achieves population stabilization.

## 2. Objectives and Methodology:

There is a major data difficulty regarding availability of annual statistics, calculations & comparisons of Crude Birth Rate (CBR), Total Fertility Rate (TFR) and Females' Mean Age at Gauna (FMAG) over time for district level study of any state and which is applied to Rajasthan also. This data problem distorts the calculations and negates the usefulness of making comparisons over time. Due to this data information problem, we use the information for different years (as per the availability of latest data, taking 2000-01 as base year) in this paper. This data problem at district level is a constraint that creates a limitation in the selection of study objectives and hypotheses. This paper attempts to identify the demographic and socio-economic variables, which are responsible for population growth in Rajasthan. The main objectives of the study are:

- ❖ To observe the characteristics of indicators of population growth in Rajasthan.
- ❖ To identify the various demographic & socio-economic variables which have causal relationship with population growth.
- ❖ To analyze the inter-relationship between the indicators of population growth and demographic & socio-economic variables.

---

[5] Government of Rajsthan (1999), "Population Policy of Rajasthan", Department of Family Welfare, Jaipur.

For achieving the above objectives, the *a priori* hypotheses are as follows:

- ❖ Positive impact of infant mortality & total fertility rate and negative impact of income equality on population growth.
- ❖ Positive impact of infant mortality and negative impact of female's age at gauna and female literacy on crude birth rate.
- ❖ Negative impact of couple protection rate, income equality, female literacy and positive impact of infant mortality on total fertility rate.
- ❖ Positive impact of female literacy & income equality on female's age at gauna.
- ❖ Positive impact of female literacy, females age at gauna and income equality on couple protection rate.

To rummage the inter-relationship between indicators of population growth and demographic & socio-economic variables, a social sector model is proposed. The model is estimated by the use of Multiple Regression Analysis (Method of Ordinary Least Squares). The general form of the Multiple Regression Equation Model is as follows:

$$Y_i = \beta_1 + \beta_2 X_{2i} + \beta_3 X_{3i} + \cdots + \beta_k X_{ki} + u_i$$
$$\text{where } i = 1, 2, 3, \ldots, n.$$

In this multiple regression equation model, $Y_i$ is dependent variable and $X_2, X_3, \ldots, X_k$ are independent explanatory variables. $\beta_1$ is the intercept, shows the average value of Y, when $X_2, X_3, \ldots, X_k$ are set equal to zero; $\beta_2, \beta_3, \ldots, \beta_k$ are partial regression/slope coefficients; $u_i$ is the stochastic disturbance term; i is the i$^{th}$ observation and n is the size of population.

The model is estimated by using cross-sectional data of all 32 districts of the state (at that time, the no. of districts was 32). In this paper, we also calculated the Mean, Standard Deviation and Coefficient of Variation of the variables. The variables used in this paper, their reference year and abbreviations/identification code are given in the Appendix I (Table 9). Firstly, we regress the dependent variables with all the variables, which have theoretical relationship and then choose the appropriate variables for multiple regressions. The dependent and independent variables for the model are as follows:

**Table 2: Functional Form of the Model**

| Dependent Variable | Independent Variables |
|---|---|
| POPGWR | CBR, TFR, FMAG, CDR, CPR, IMR, CIMM, MRANC, PWRSAP, PWETVR, MIPLP, BPGH, PCEMPH, LIT, LIT$_m$, LIT$_f$, PCEEE, PCNDDP, PPBPL, ROADSK, PHDW, PCEWS |
| CBR | POPGWR, FMR, FMR$_{(0-6)}$, PURPOP, FMAG, CPR, IMR, PWETVR, PCEMPH, PCEFW, LIT, LIT$_m$, LIT$_f$, PCEEE, PCNDDP, PPBPL |
| TFR | PURPOP, FMAG, CPR, CDR, IMR, MRANC, PWRSAP, PWETVR, PCEMPH, LIT, LIT$_m$, LIT$_f$, PCEEE, PCNDDP, PPBPL, PCESCS |
| FMAG | PURPOP, PWETVR, LIT, LIT$_m$, LIT$_f$, PSER, PSER$_m$ PSER$_f$, DORPS, DORPS$_m$, DORPS$_f$, PCEEE, PCNDDP, PPBPL |

| Dependent Variable | Independent Variables |
|---|---|
| CPR | PURPOP, FMAG, IMR, PWETVR, MIPLP, PCEMPH, PCEFW, LIT, $LIT_m$, $LIT_f$, PCEEE, PCNDDP, PPBPL, IDI, PCESCS |

In this paper, we have taken 32 variables (appendix-I). All the 32 variables are relating to Population; Fertility, Reproductive Health and Mortality; Public Health and Health Infrastructure; Education and Educational Infrastructure; and Economic Growth and Infrastructure. Data used in this paper have taken from website of Census Department, State Human Development Report (Rajasthan), Various Administrative Reports of Medical, Health & Family Welfare Department, Government of Rajasthan and Plan Documents of Planning Department, Government of Rajasthan.

### 3. Multivariate Analysis

### 3.1 Mean, Standard Deviation & Coefficient of Variation

Mean, standard deviation and coefficient of variation of all the 32 variables for all 32 districts along with the figures of all Rajasthan are at appendix I (table 9). The Mean, measures the average value of the variables for all 32 districts. The Standard Deviation, measures the absolute variation in the mean and the Coefficient of Variation, measures the percentage variation in mean. The variables are divided in to five categories according to the range of Coefficient of Variation for the analysis of Standard Deviation and Coefficient of Variation.

**Table 3: Range-wise Variables according to the Coefficient of Variation**

| Range | Variables |
|---|---|
| Less than 25% | POPGWR (19.92), FMR (5.28), $FMR_{(0-6)}$(3.26), CBR (7.02), TFR (10.20), FMAG (3.66), CPR (14.73), CDR (10.44), IMR (20.60), MIPLP (18.59), LIT (12.64), $LIT_m$ (8.31), $LIT_f$ (21.19), PSER (9.39), $PSER_m$ (11.17), $PSER_f$ (14.27), DORPS (12.94), $DORPS_m$ (12.89), $DORPS_f$ (17.01), PCNDDP (24.34), PHDW (21.33) |
| 25% to 50% | MRANC (38.95), BPGH (30.61), PCEFW (38.55), PCEEE (44.54), PPBPL (46.88), PCESCS (48.75) |
| 50% to 75% | PURPOP (53.79), PWETVR (66.94), PCEMPH (58.63), IDI (55.05) |
| 75% to 100% | - |
| More than 100% | PCEWS (158.62) |

Table 3 shows that variability is higher in the variables of public health & health infrastructure and economic growth & infrastructure head. There is need to reduce disparities on this front.

### 3.2 Regression analysis

To rummage the interrelationship between indicators of population growth and various demographic and socio-economic variables, we regress the dependent variable with the independent variables individually (independent variables are those variable which have causal relationship with dependent variable in theoretical and behavioral terms) and then pick the most influential variables and regress with the help of step-wise method and get best fitted multiple regression equation of them. Some variables with insignificant coefficients have also been kept in the model because theoretically their importance has been proved. Figures below the coefficients are 't' values. Significance of variables with the level of significance is denoted as follows:

    \*        Significant at 1% level of significance
    \*\*      Significant at 2% level of significance
    \*\*\*     Significant at 5% level of significance
    \*\*\*\*    Significant at 10% level of significance

Efforts have been made to avoid the problem of multicollinearity (as it presents commonly in the analysis of cross-sectional data) but at some places, it is difficult to avoid it.

### 3.2.1 Population Growth (Decennial)

Population Growth (POPGWR) is regressed with different variables such as CBR, TFR, FMAG, CDR, CPR, IMR, CIMM, MRANC, PWRSAP, PWETVR, MIPLP, BPGH, PCEMPH, LIT, $LIT_m$, $LIT_f$, PCNDDP, PPBPL, ROADSK, PHDW, PCEWS.

**Table 4: Regression Equations of Population Growth (Decennial)**

| S.No. | Intercept | | Coefficient | | $R^2$ | d. f. |
|---|---|---|---|---|---|---|
| 1. | 10.0934 | + | 0.5643 | CBR | 0.0514 | 31 |
| | | | 1.2745 | | | |
| 2. | 8.1549 | + | 4.1092 | TFR*** | 0.1325 | 31 |
| | | | 2.1408 | | | |
| 3. | 49.0047 | − | 1.1555 | FMAG | 0.0156 | 31 |
| | | | 0.6903 | | | |
| 4. | 50.0223 | − | 0.5750 | CPR* | 0.3245 | 31 |
| | | | 3.7963 | | | |
| 5. | 46.2904 | − | 2.0212 | CDR**** | 0.1119 | 31 |
| | | | 1.9442 | | | |
| 6. | 36.6587 | − | 0.0979 | IMR**** | 0.0946 | 31 |
| | | | 1.7709 | | | |
| 7. | 34.4649 | − | 0.1670 | CIMM*** | 0.1413 | 31 |
| | | | 2.2217 | | | |
| 8. | 278632 | + | 0.0061 | MRANC | 0.0007 | 31 |
| | | | 0.1473 | | | |
| 9. | 13.8825 | + | 0.1462 | PWRSAP | 0.0497 | 31 |
| | | | 1.2532 | | | |
| 10. | 30.8793 | − | 0.1961 | PWETVR**** | 0.097 | 31 |
| | | | 1.8019 | | | |
| 11. | 24.9574 | + | 0.1171 | MIPLP | 0.0118 | 31 |
| | | | 0.5995 | | | |
| 12. | 21.9702 | + | 0.0768 | BPGH**** | 0.1167 | 31 |
| | | | 1.9906 | | | |
| 13. | 29.1035 | − | 0.0453 | PCEMPH | 0.0079 | 31 |
| | | | 0.4881 | | | |
| 14. | 32.8303 | − | 0.0768 | LIT | 0.0106 | 31 |
| | | | 0.5664 | | | |
| 15. | 40.5194 | − | 0.1629 | $LIT_m$ | 0.0328 | 31 |
| | | | 1.0084 | | | |
| 16. | 31.2097 | − | 0.0696 | $LIT_f$ | 0.0124 | 31 |
| | | | 0.6137 | | | |
| 17. | 26.7674 | + | 0.0346 | PCEEE | 0.0138 | 31 |
| | | | 0.6478 | | | |
| 18. | 32.6477 | − | 0.0003 | PCNDDP | 0.0361 | 31 |
| | | | 1.0604 | | | |
| 19. | 27.1559 | + | 0.0285 | PPBPL | 0.0057 | 31 |
| | | | 0.4138 | | | |
| 20. | 35.2376 | − | 0.2308 | ROADSK*** | 0.1539 | 31 |
| | | | 2.3363 | | | |
| 21. | 31.0646 | − | 0.0465 | PHDW | 0.0114 | 31 |
| | | | 0.5872 | | | |
| 22. | 28.6427 | − | 0.0134 | PCEWS | 0.0122 | 31 |
| | | | 0.6075 | | | |

Fit of the equations is with the expected signs. TFR, CPR, CDR, IMR, CIMM, PWETVR, BPGH and ROADSK have significant coefficients. PCEEE appears with opposite sign as of expected sign. In the step-wise regression, PPBPL is found more relevant in spite of PCNDDP for multiple regression.

$$POPGWR = 12.5485 + 5.6405 \text{ TFR*} - 0.1477 \text{ IMR*} + 0.0246 \text{ PPBPL}$$
$$(3.0425) \quad\quad (2.7565) \quad\quad (0.4075)$$
$$R^2 = 0.3196 \quad\quad d.f. = 29$$

In the multiple regression analysis the coefficients of TFR and IMR are significant at 1% level of significance. This indicates that TFR influences POPGWR positively. IMR shows negative influence to POPGWR in mathematical/statistical terms but in actual terms this leads to birth to more children due to less survival. The variable PPBPL does not affect POPGWR significantly.

### 3.2.2 Crude Birth Rate

Crude Birth Rate (CBR) is regressed with different variables such as POPGWR, FMR, $FMR_{(0-6)}$, PURPOP, FMAG, CPR, IMR, PWETVR, PCEMPH, PCEFW, LIT, $LIT_m$, $LIT_f$, PCEEE, PCNDDP, PPBPL.

**Table 5: Regression Equations of Crude Birth Rate**

| S.No. | Intercept | | Coefficient | | $R^2$ | d. f. |
|---|---|---|---|---|---|---|
| 1. | 29.6053 | + | 0.0910 | POPGWR | 0.0514 | 31 |
| | | | 1.2745 | | | |
| 2. | 39.4165 | - | 0.0079 | FMR | 0.0286 | 31 |
| | | | 0.9391 | | | |
| 3. | 25.6612 | - | 0.0072 | $FMR_{(0-6)}$ | 0.0088 | 31 |
| | | | 0.5163 | | | |
| 4. | 32.1109 | + | 0.0032 | PURPOP | 0.0002 | 31 |
| | | | 0.0856 | | | |
| 5. | 50.5045 | - | 1.1004 | FMAG**** | 0.0879 | 31 |
| | | | 1.7004 | | | |
| 6. | 38.4234 | - | 0.1650 | CPR*** | 0.1657 | 31 |
| | | | 2.4406 | | | |
| 7. | 30.0598 | + | 0.0247 | IMR | 0.0372 | 31 |
| | | | 1.0766 | | | |
| 8. | 32.2562 | - | 0.0059 | PWETVR | 0.0006 | 31 |
| | | | 0.1291 | | | |
| 9. | 31.1170 | + | 0.0563 | PCEMPH | 0.0755 | 31 |
| | | | 1.5657 | | | |
| 10. | 32.3631 | - | 0.1645 | PCEFW | 0.0010 | 31 |
| | | | 0.1744 | | | |
| 11. | 32.4349 | - | 0.0043 | LIT | 0.0002 | 31 |
| | | | 0.0792 | | | |
| 12. | 30.2455 | + | 0.0256 | $LIT_m$ | 0.0050 | 31 |
| | | | 0.3897 | | | |
| 13. | 32.9059 | - | 0.0172 | $LIT_f$ | 0.0047 | 31 |
| | | | 0.3753 | | | |
| 14. | 32.2458 | - | 0.0016 | PCEEE | 0.0002 | 31 |
| | | | 0.0748 | | | |
| 15. | 34.6147 | - | 0.0002 | PCNDDP | 0.0689 | 31 |
| | | | 1.4901 | | | |
| 16. | 31.1572 | + | 0.0321 | PPBPL | 0.0447 | 31 |
| | | | 1.1847 | | | |

FMAG and CPR have significant coefficients. PURPOP, PCEMPH and $LIT_m$ are with opposite signs as of expected signs.

$$CBR = 50.2161 - 1.0819 \ FMAG^{****} + 0.0114 \ IMR - 0.0234 \ LIT_f$$
$$(1.7123) \qquad (0.4421) \qquad (0.4701)$$
$$R^2 = 0.1099 \qquad d.f. = 29$$

Fit of the multiple regression equation is with the expected signs Coefficient of FMAG is significant at 10% level of significance. This indicates that FMAG influences CBR negatively. The coefficients of IMR and $LIT_f$ are insignificant but included due to their importance in the determination of CBR.

### 3.2.3 Total Fertility Rate

Total Fertility Rate (TFR) is regressed with different variables such as PURPOP, FMAG, CPR, CDR, IMR, MRANC, PWRSAP, PWETVR, PCEMPH, LIT, $LIT_m$, $LIT_f$, PCEEE, PCNDDP, PPBPL, PCESCS.

**Table 6: Regression Equations of Total Fertility Rate**

| S.No. | Intercept | | Coefficient | | $R^2$ | d. f. |
|---|---|---|---|---|---|---|
| 1. | 4.9033 | - | 0.0006 | PURPOP | 0.0002 | 31 |
| | | | 0.0744 | | | |
| 2. | 9.2697 | - | 0.2629 | FMAG**** | 1.1031 | 31 |
| | | | 1.8573 | | | |
| 3. | 6.5736 | - | 0.0445 | CPR* | 0.2471 | 31 |
| | | | 3.1379 | | | |
| 4. | 3.6639 | + | 0.1375 | CDR | 0.0659 | 31 |
| | | | 1.7123 | | | |
| 5. | 4.2069 | + | 0.0079 | IMR | 0.0797 | 31 |
| | | | 1.6123 | | | |
| 6. | 5.2989 | - | 0.0064 | MRANC**** | 0.1017 | 31 |
| | | | 1.8431 | | | |
| 7. | 5.0179 | - | 0.0033 | PWRSAP | 0.0032 | 31 |
| | | | 0.3124 | | | |
| 8. | 51792 | - | 0.0073 | PWETVR | 0.0175 | 31 |
| | | | 0.7304 | | | |
| 9. | 4.9694 | - | 0.0011 | PCEMPH | 0.0006 | 31 |
| | | | 0.1367 | | | |
| 10. | 4.9951 | - | 0.0033 | LIT | 0.0025 | 31 |
| | | | 0.2719 | | | |
| 11. | 5.4818 | - | 0.0139 | $LIT_m$ | 0.0305 | 31 |
| | | | 0.9720 | | | |
| 12. | 5.0591 | - | 0.0039 | $LIT_f$ | 0.0051 | 31 |
| | | | 0.3932 | | | |
| 13. | 4.9577 | - | 0.0016 | PCEEE | 0.0036 | 31 |
| | | | 0.3288 | | | |
| 14. | 5.6749 | - | 0.00006 | PCNDDP*** | 0.1465 | 31 |
| | | | 2.2693 | | | |
| 15. | 4.9662 | + | 0.0023 | PPBPL | 0.0051 | 31 |
| | | | 0.3906 | | | |
| 16. | 5.1543 | - | 0.0014 | PCESCS | 0.0664 | 31 |
| | | | 1.4618 | | | |

FMAG, CPR, MRANC and PCNDDP are with significant coefficients. All the variables show the expected signs.

$$TFR = 5.9697 - 0.0412\ CPR^* + 0.0104\ IMR^{***}$$
$$(2.9361) \quad\quad (2.3704)$$
$$- 0.0031\ LIT_f - 0.00004\ PCNDDP^{****}$$
$$(0.3446) \quad\quad (1.8641)$$
$$R^2 = 0.4305 \quad\quad d.f. = 28$$

Coefficient of CPR is significant at 1% level of significance, IMR at 2% and PCNDDP at 10%. This indicates that CPR & PCNDDP influence TFR positively and IMR influences TFR negatively. $LIT_f$ appears with insignificant coefficient but it has major influential role in the determination of TFR.

### 3.2.4 Females' Mean Age at Gauna

Females' Mean Age at Gauna (FMAG) is regressed with different variables such as PURPOP, PWETVR, LIT, $LIT_m$, $LIT_f$, PSER, $PSER_m$ $PSER_f$, DORPS, $DORPS_m$, $DORPS_f$, PCEEE, PCNDDP, PPBPL.

**Table 7: Regression Equations of Females' Mean Age at Gauna**

| S.No. | Intercept | | Coefficient | | $R^2$ | d.f. |
|---|---|---|---|---|---|---|
| 1. | 16.1794 | + | 0.0060 | PURPOP | 0.0118 | 31 |
| | | | 0.5994 | | | |
| 2. | 15.0223 | + | 0.0273 | PWETVR*** | 0.1619 | 31 |
| | | | 2.4070 | | | |
| 3. | 15.7873 | + | 0.0190 | LIT | 0.051 | 31 |
| | | | 1.3231 | | | |
| 4. | 14.9522 | + | 0.0305 | $LIT_m$**** | 0.0981 | 31 |
| | | | 1.8061 | | | |
| 5. | 151910 | + | 0.0126 | $LIT_f$ | 0.0346 | 31 |
| | | | 1.0371 | | | |
| 6. | 16.0412 | + | 0.0160 | PSER | 0.0456 | 31 |
| | | | 1.1974 | | | |
| 7. | 16.9408 | − | 0.0028 | $PSER_m$ | 0.0027 | 31 |
| | | | 0.2861 | | | |
| 8. | 15.8440 | + | 0.0166 | $PSER_f$ | 0.0775 | 31 |
| | | | 1.5871 | | | |
| 9. | 17.3288 | − | 0.0225 | DORPS | 0.0796 | 31 |
| | | | 1.6104 | | | |
| 10. | 18.1252 | − | 0.0268 | $DORPS_m$**** | 0.1050 | 31 |
| | | | 1.8765 | | | |
| 11. | 17.2214 | − | 0.0099 | $DORPS_f$ | 0.0147 | 31 |
| | | | 0.6701 | | | |
| 12. | 16.7143 | + | 0.0014 | PCEEE | 0.0018 | 31 |
| | | | 0.2328 | | | |
| 13. | 16.4417 | + | 0.00002 | PCNDDP | 0.0074 | 31 |
| | | | 0.4713 | | | |
| 14. | 16.6888 | − | 0.0010 | PPBPL | 0.0006 | 31 |
| | | | 0.1374 | | | |

PWETVR, $LIT_m$ and $DORPS_m$ are with significant coefficients. Except $PSER_m$, coefficients of all are with expected Signs.

$$FMAG = 13.7224 + 0.0279\ PWETVR^{***} + 0.0039\ LIT_f^{****} + 0.00003\ PCNDDP$$
$$(2.1774) \qquad\qquad (1.8126) \qquad\qquad (1.0034)$$
$$R^2 = 0.1912 \qquad d.f. = 29$$

All the variables are with expected signs. Coefficient of PWETVR is significant at 5% level of significance & coefficient of $LIT_f$ is significant at 10% level of significance. This indicates that PWETVR & $LIT_f$ influence FMAG positively. Coefficient of PCNDDP is insignificant means the variable PCNDDP does not affect FMAG significantly.

### 3.2.5 Couple Protection Rate

Couple Protection Rate (CPR) is regressed on different variables such as PURPOP, FMAG, IMR, PWETVR, MIPLP, PCEMPH, PCEFW, LIT, $LIT_m$, $LIT_f$, PCEEE, PCNDDP, PPBPL, IDI, PCESCS.

**Table 8: Regression Equations of Couple Protection Rate**

| S.No. | Intercept | | Coefficient | | $R^2$ | d.f. |
|---|---|---|---|---|---|---|
| 1. | 40.2031 | − | 0.1133 | PURPOP | 0.0511 | 31 |
| | | | 1.2713 | | | |
| 2. | 18.8647 | + | 1.1404 | FMAG | 0.0155 | 31 |
| | | | 0.6876 | | | |
| 3. | 34.1252 | + | 0.0435 | IMR | 0.0190 | 31 |
| | | | 0.7628 | | | |
| 4. | 35.2834 | + | 0.1922 | PWETVR**** | 0.0956 | 31 |
| | | | 1.7811 | | | |
| 5. | 35.3518 | + | 0.0892 | MIPLP | 0.0069 | 31 |

| S.No. | Intercept | | Coefficient | | $R^2$ | d. f. |
|---|---|---|---|---|---|---|
| | | | 0.4596 | | | |
| 6. | 36.7349 | + | 0.0597 | PCEMPH | 0.0139 | 31 |
| | | | 0.6522 | | | |
| 7. | 33.9617 | + | 3.4369 | PCEFW | 0.0726 | 31 |
| | | | 1.5325 | | | |
| 8. | 35.5692 | + | 0.1294 | LIT | 0.0306 | 31 |
| | | | 0.9726 | | | |
| 9. | 29.9513 | + | 0.1606 | $LIT_m$ | 0.0325 | 31 |
| | | | 1.0032 | | | |
| 10. | 36.5518 | + | 0.0869 | $LIT_f$ | 0.0197 | 31 |
| | | | 0.7762 | | | |
| 11. | 39.5813 | + | 0.0402 | PCEEE**** | 0.0189 | 31 |
| | | | 1.7607 | | | |
| 12. | 31.0288 | + | 0.0005 | PCNDDP**** | 0.0889 | 31 |
| | | | 1.7108 | | | |
| 13. | 39.0037 | - | 0.1215 | PPBPL**** | 0.1051 | 31 |
| | | | 1.8769 | | | |
| 14. | 38.5776 | + | 0.0077 | IDI | 0.0050 | 31 |
| | | | 0.3894 | | | |
| 15. | 36.6705 | + | 0.0093 | PCESCS | 0.0250 | 31 |
| | | | 0.8785 | | | |

PWETVR, PCEEE, PCNDDP and PPBPL are with significant coefficients and expected signs. Sign of coefficient of PURPOP is opposite of the expected.

$$CPR = 20.6541 + 0.4813 \text{ FMAG} + 0.1388 \text{ LIT}_f{}^{****} + 0.0006 \text{ PCNDDP}^{****}$$
$$(0.2922) \qquad (1.8065) \qquad (1.9266)$$
$$R^2 = 0.1433 \qquad d.f. = 29$$

All the variables are with expected signs of coefficients. Coefficients of $LIT_f$ and PCNDDP are significant at 10% level of significance. This indicates that $LIT_f$ and PCNDDP influence CPR positively. Coefficient of FMAG is insignificant means the variable FMAG does not affect CPR significantly.

## 4.    Conclusion

The model is fit good with the expected signs. Estimated equations confirm the *a priori* hypotheses of positive impact of infant mortality & total fertility rate and negative impact of income equality on population growth; positive impact of female literacy & income equality on female's age at gauna; positive impact of infant mortality and negative impact of female's age at gauna and female literacy on crude birth rate; negative impact of couple protection rate, income equality, female literacy and positive impact of infant mortality on total fertility rate, positive impact of female literacy, females age at gauna and income equality on couple protection rate. Literacy, especially female literacy and per-capita income appeared as most influential variables to attack the poor status of socio-economic & demographic variables. There is need to emphasize on the improvement of these two variables.

Rapid population growth retards the economic, social and human development. Enhancement of women's status and autonomy has been conclusively established to have a direct bearing on fertility and mortality decline, which indirectly affects the population growth. More specifically, inter-relationships between women's characteristics and access to resources are the mechanisms through which human fertility is determined. Education is highly correlated with age at the marriage of the females and thus helps in the reduction of the reproductive life, on an average, and helps in the conscious efforts to limit the family size. The early marriage of the daughter in rural areas is an expected rational behavior, as long as there is mass illiteracy and poverty. The age at marriage for females cannot be raised by mere, legislation unless the socio-economic conditions of the rural people is improved and better educational facilities and occupational alternatives for the teenage girls are provided near their homes.

Reproductive and public health have their importance in determination of population stabilization. National Rural Health Mission (NRHM) and Rajasthan Health System Development Project (RHSDP) are ongoing programs which can improve the situation. There is need of effective monitoring of activities under these programs. Effective implementation of family welfare program will create opportunities for better education and improvement in nutritional status of family through check on population growth, which will turn in better health of mother and child and there will be less infant and maternal mortality.

**References**

- Government of Rajasthan (2005), "District-wise Performance of Family Welfare Programme-2004", Directorate of Family Welfare, Jaipur.
- Government of Rajasthan, "Various Plan Documents", Planning Department, Jaipur.
- Government of Rajsthan (1999), "Population Policy of Rajasthan", Department of Family Welfare, Jaipur.
- Kulkarni, Sumati and Minja Kim Choe (1997), "State-level Variations in Wanted and Unwanted Fertility Provide a Guide for India's Family Planning Programmes", NFHS Bulletion, IIPS, Mumbai.
- Mittal, Alka (2004), "Billion Plus Population: Challenges Ahead", Paper submitted to Academic Staff College, University of Rajasthan, Jaipur during 57th Orientation Course.
- Mohanty, Sanjay K. and Moulasha K. (1996), "Women's Status, Proximate Determinants and Fertility Behaviour in Rajasthan", Paper Presented at National Seminar on Population and Development in Rajasthan at HCM-RIPA, Jaipur.
- Murthy, M.N. (1996), "Reasons for Low Contraceptive Use in Rajasthan", Paper Presented at National Seminar on Population and Development in Rajasthan at HCM-RIPA, Jaipur.
- Radhakrishan, S., S. Sureender and R. Acharya (1996), "Child Marriage: Determinants and Attitudes Examined in Rajasthan", Paper Presented at National Seminar on Population and Development in Rajasthan at HCM-RIPA, Jaipur.
- Ramesh, B.M., S.C. Gulati and Robert D. Retherford (1996), "Contraceptive Use in India", NFHS Subject Report, IIPS, Mumbai.
- Retherford, Robert D. and Vinod Mishra (1997), "Media Exposure Increases Contraceptive Use", NFHS Bulletin, IIPS, Mumbai.
- Retherford, Robert D., M.M. Gandotra, Arvind Pandey, Norman Y. Luther, and Vinod K. Mishra (1998), "Fertility in India", NFHS Subject Report, IIPS, Mumbai.
- Retherford, Robert D., P.S. Nair, Griffith Feeney and Vinod K. Mishra (1999), "Factors Affecting Source of Family Planning Services in India", NFHS Subject Report, IIPS, Mumbai.
- Roy, T.K., R. Mutharayappa, Minja Kim choe and Fred Arnold (1997), "Son Preference and its Effect on Fertility in India", NFHS Subject Report, IIPS, Mumbai.
- Shariff, Abusaleh (1996), "Poverty and Fertility Differentials in Indian States: New Evidence from Cross-Sectional Data", Margin, October-December, Vol. 29, No.1, pp. 49-67.
- Sinha, Narain and Assakaf Ali (1999), "Econometric Analysis of Socio-Economic Determinants of Fertility: A Case Study of Yemen", Paper Presented at the Conference of the India Econometric Society, Jaipur.
- Society for International Development (1999), "Human Development Report: Rajasthan", Rajasthan Chapter, Jaipur.
- Visaraia, Pravin and Leela Visaria (1995), "India's Population in Transition", Population Bulletin, 50(3), Population Reference Bureau, Washington, D.C.
- website www.censusindia.net

**Table 9: All Rajasthan Figures, Mean, Standard Deviation & Coefficient of Variation of Variables**

| S. No. | Variable & Year | Code | Unit | All Rajasthan | Mean | S. D. | CoV |
|---|---|---|---|---|---|---|---|
| 1. | Population Growth (Decennial) 1991-2001 | POPGWR | Per cent | 28.33 | 28.25 | 5.63 | 19.92 |
| 2. | Female-Male Ratio 2001 | FMR | Nos. | 921 | 922.03 | 48.65 | 5.28 |
| 3. | Female-Male Ratio (0-6 years) 2001 | $FMR_{(0-6)}$ | Nos. | 909 | 909.00 | 29.59 | 3.26 |
| 4. | Percentage of Urban Population to Total Population 2001 | PURPOP | Per cent | 23.38 | 20.69 | 11.13 | 53.79 |
| 5. | Crude Birth Rate 1997 | CBR | Per '000 | 32.90 | 32.18 | 2.26 | 7.02 |
| 6. | Total Fertility Rate 1997 | TFR | Nos. | 4.9 | 4.89 | 0.50 | 10.20 |
| 7. | Females Mean Age at Gauna 1996-97 | FMAG | Years | 17.7 | 16.66 | 0.61 | 3.66 |
| 8. | Couple Protection Rate 2001 | CPR | Per cent | 37.00 | 37.86 | 5.58 | 14.73 |
| 9. | Crude Death Rate 1997 | CDR | Per '000 | 8.9 | 8.93 | 0.93 | 10.44 |
| 10. | Infant Mortality Rate 1997 | IMR | Per '000 | 87 | 85.81 | 17.67 | 20.60 |
| 11. | Percentage of Mothers Receiving Total Ante-Natal Care 1996-97 | MRANC | Per cent | 72.3 | 63.38 | 24.69 | 38.95 |
| 12. | Percentage of Women having Exposure to TV & Radio 1996-97 | PWETVR | Per cent | 13.1 | 13.40 | 8.97 | 66.94 |
| 13. | Medical Institutions Per-Lakh of Population 1997-98 | MIPLP | Nos. | 27 | 28.13 | 5.23 | 18.59 |
| 14. | Beds Per-Lakh Population in Govt. Hospitals 1997-98 | BPGH | Nos. | 85 | 81.81 | 25.04 | 30.61 |
| 15. | Per-Capita Expenditure on Medical & Public Health 2000-01 | PCEMPH | ` | 19.00 | 18.82 | 11.04 | 58.63 |
| 16. | Per-Capita Expenditure on Family Welfare 2000-01 | PCEFW | ` | 0.97 | 1.13 | 0.44 | 38.55 |
| 17. | Literacy Rate 2001 | LIT | Per cent | 60.41 | 59.58 | 7.53 | 12.64 |
| 18. | Literacy Rate (Male) 2001 | $LIT_m$ | Per cent | 75.70 | 75.31 | 6.26 | 8.31 |
| 19. | Literacy Rate (Female) 2001 | $LIT_f$ | Per cent | 43.85 | 42.51 | 9.01 | 21.19 |
| 20. | Primary School Enrolment Ratio 1997-98 | PSER | Per cent | 86.50 | 86.75 | 8.15 | 9.39 |
| 21. | Primary School Enrolment Ratio (Male) 1997-98 | $PSER_m$ | Per cent | 99.78 | 100.51 | 11.22 | 11.17 |
| 22. | Primary School Enrolment Ratio (Female) 1997-98 | $PSER_f$ | Per cent | 71.91 | 71.65 | 10.22 | 14.27 |
| 23. | Drop-Out Rates at Primary Level 1996-97 | DORPS | Per cent | 56.60 | 59.13 | 7.65 | 12.94 |
| 24. | Drop-Out Rates at Primary Level (Male) 1996-97 | $DORPS_m$ | Per cent | 54.72 | 57.07 | 7.36 | 12.89 |
| 25. | Drop-Out Rates at Primary Level (Female) 1996-97 | $DORPS_f$ | Per cent | 56.96 | 62.68 | 10.66 | 17.01 |
| 26. | Per-Capita Expenditure on Elementary Education 2000-01 | PCEPEE | ` | 47.00 | 42.86 | 19.09 | 44.54 |
| 27. | Per-Capita Net District Domestic Product 1999-2000 | PCNDDP | ` | 12752 | 12831.88 | 3122.80 | 24.34 |
| 28. | Population Below Poverty Line 1999-2000 | PPBPL | Per cent | 30.99 | 31.74 | 14.88 | 46.88 |
| 29. | Infrastructure Development Index 1994-95 | IDI | Nos. | 100.00 | 93.46 | 51.45 | 55.05 |
| 30. | Percentage of Villages with Safe Drinking Water 1998-99 | PHDW | Per cent | 64.30 | 60.54 | 12.91 | 21.33 |
| 31. | Per-Capita Expenditure on Social & Community Services 2000-01 | PCESCS | ` | 245.62 | 194.69 | 94.92 | 48.75 |
| 32. | Per-Capita Expenditure on Water Supply 2000-01 | PCEWS | ` | 39.95 | 29.19 | 46.30 | 158.62 |